\documentclass{llncs}

\usepackage{t1enc}
\usepackage[latin2]{inputenc}
\usepackage[english]{babel}
\usepackage{url}
\usepackage{amsmath}
\usepackage[mathscr]{eucal}
\usepackage{subfig}
\usepackage{bm}
\usepackage{ntheorem}
\usepackage{amssymb}
\usepackage{graphicx}

\newcommand{\R}{\ensuremath{\mathbb{R}}}
\newcommand{\C}{\ensuremath{\mathbb{C}}}

\newtheorem*{theorem2}{Theorem}
\newtheorem*{proof2}{Proof}

\renewcommand{\b}{\mathbf}
\title{Independent Process Analysis without A Priori Dimensional Information}
\author{Barnab\'as P\'oczos, Zolt\'an Szab\'o, Melinda Kiszlinger, and Andr\'as L{\H o}rincz\thanks{Corresponding author}}
\institute{Department of Information Systems \\
E\"{o}tv\"{o}s Lor\'{a}nd University, Budapest, Hungary \\
\email{\{pbarn,szzoli\}@cs.elte.hu, \{kmelinda,andras.lorincz\}@elte.hu}}
\begin{document}
\date{}
\maketitle

\begin{abstract}
Recently, several algorithms have been proposed for independent subspace analysis
where hidden variables are i.i.d. processes. We show that these methods can be
extended to certain AR, MA, ARMA and ARIMA tasks. Central to our paper is that we
introduce a cascade of algorithms, which aims to solve these tasks without previous
knowledge about the number and the dimensions of the hidden processes. Our claim is
supported by numerical simulations. As a particular application, we search for
subspaces of facial components.
\end{abstract}

\section{Introduction}
\label{intro} Independent Subspace Analysis (ISA), also known as Multidimensional
Independent Component Analysis \cite{cardoso98multidimensional}, is a generalization
of Independent Component Analysis (ICA). ISA assumes that certain sources depend on
each other, but the dependent groups of sources are still independent of each other,
i.e., the independent groups are multidimensional. The ISA task has been subject of
extensive research
\cite{cardoso98multidimensional,hyvarinen00emergence,vollgraf01multi,strogbauer04least,bach03beyond,theis05blind,hyvarinnen06fastisa,nolte06identifying,poczos05geodesic}.
In this case, one assumes that the hidden sources are independent and identically
distributed (i.i.d.) in time. Temporal independence is, however, a gross
oversimplification of real sources including acoustic or biomedical data. One may
try to overcome this problem, by assuming that hidden processes are, e.g.,
autoregressive (AR) processes. Then we arrive to the AR Independent Process Analysis
(AR-IPA) task \cite{hyvarinen98stochastic,poczos05independent}. Another method to
weaken the i.i.d. assumption is to assume moving averaging (MA). This direction is
called Blind Source Deconvolution (BSD) \cite{choi05blind}, in this case the
observation is a temporal mixture of the i.i.d. components.

The AR and MA models can be generalized and one may assume ARMA sources instead of
i.i.d. ones. As an additional step, the method can be extended to non-stationary
integrated ARMA (ARIMA) processes, which are important, e.g., for modelling economic
processes \cite{mills90time}.

In this paper, we formulate the AR-, MA-, ARMA-, ARIMA-IPA generalization of the ISA tasks, when (i) one allows for
multidimensional hidden components and (ii) the dimensions of the hidden processes are not known. We show that in the
undercomplete case, when the number of `sensors' is larger than the number of `sources', these tasks can be reduced to
the ISA task.

\section{Independent Subspace Analysis \label{s:ISA}} The ISA task can be formalized as follows:
\begin{align}
    \b{x}(t)=\b{A}\b{e}(t), \text{ where }
    \b{e}(t)=\left[\b{e}^1(t);\ldots;\b{e}^M(t)\right]\in\R^{D_e} \label{eq:e_concat}
\end{align}
and $\b{e}(t)$ is a vector concatenated of components
$\b{e}^m(t)\in\R^{d_e^m}$. The total dimension of the components
is $D_e=\sum_{m=1}^M d_e^m$. We assume that for a given $m$,
$\b{e}^m(t)$ is i.i.d. in time $t$, and sources $\b{e}^m$ jointly
independent, i.e., $I(\b{e}^1,\ldots,\b{e}^M)=0$, where $I(.)$
denotes the mutual information (MI) of the arguments. The
dimension of observation $\b{x}$ is $D_x$. Assume that $D_x>D_e$,
and $\b{A} \in \R^{D_x \times D_e}$ is of full column rank. Under
these conditions, one may assume without any loss of generality
that both the observed ($\b{x}$) and the hidden ($\b{e}$) signals
are white. For example, one may apply Principal Component Analysis
(PCA) as a preprocessing stage. Then the ambiguities of the ISA
task are as follows \cite{theis04uniqueness}: Sources can be
determined up to permutation and up to orthogonal transformations
within the subspaces.

\subsection{The ISA Separation Theorem}

We are to uncover the independent subspaces. Our task is to find a
matrix $\b{W} \in \R^{D_e \times D_x}$ such that
$\b{y}(t)=\b{W}\b{x}(t)$,
$\b{y}(t)=\left[\b{y}^1(t);\ldots;\b{y}^M(t)\right]$,
$\b{y}^m=[y^m_1;\ldots;y^m_{d^m_e}] \in \R^{d_e^m}$,
$(m=1,\ldots,M)$ with the condition that components $\b{y}^m$ are
independent. Here, (i) $y^m_i$ denotes the $i^{th}$ coordinate of
the $m^{th}$ estimated subspace, and (ii) $\b{W}$ can be chosen to
be orthogonal because of the whitening assumption. This task can
be solved by means of cost function that aims to minimize the
mutual information between components:
\begin{align}
J_1(\b{W})&\doteq I(\b{y}^1,\ldots,\b{y}^M). \label{e:J_1}
\end{align}
One can rewrite $J_1(\b{W})$ as follows:
\begin{align}
J_2(\b{W})&\doteq I(y_1^1,\ldots,y_{d^M_e}^M)-\sum_{m=1}^M
I(y_1^m,\ldots,y_{d^m_e}^m). \label{e:J_2}
\end{align}
The first term of the r.h.s. is the ICA cost function; it aims to
minimize mutual information for all coordinates. The other term is
a kind of \emph{anti-ICA} term; it aims to maximize mutual
information within the subspaces. One may try to apply a
heuristics and to optimize \eqref{e:J_2} in order: (1) Start by
any 'infomax' ICA algorithm and minimize the first term of the
r.h.s. in \eqref{e:J_2}. (2) Apply only permutations to the
coordinates such that they optimize the second term. In this
second step coordinates are not changed, but \eqref{e:J_2} may
decrease further. Surprisingly, this heuristics leads to the
global minimum of \eqref{e:J_1} in many cases. In other words, in
many cases, ICA that minimizes the first term of the r.h.s. of
\eqref{e:J_2} solves the ISA task apart from the grouping of the
coordinates into subspaces. This feature was observed by Cardoso,
first \cite{cardoso98multidimensional}. The extent of this feature
is still an open issue. Nonetheless, we call it `\emph{Separation
Theorem}', because for elliptically symmetric sources and for some
other distribution types one can prove that it is rigorously true
\cite{szabo06separation}. (See also, the result concerning local
minimum points \cite{theis06towards}). Although there is no proof
for general sources as of yet, a number of algorithms applies this
heuristics with success
\cite{cardoso98multidimensional,strogbauer04least,theis06towards,bach03finding,szabo06cross,meraim04algorithms}.
\subsection{ISA with Unknown Components}

Another issue concerns the computation of the second term of \eqref{e:J_2}. If the
$d_e^m$ dimensions of subspaces $\b{e}^m$ are known then one might rely on
multi-dimensional entropy estimations \cite{poczos05geodesic}, but these are
computationally expensive. Other methods deal with implicit or explicit pair-wise
dependency estimations \cite{bach03finding,theis06towards}. Interestingly, if the
observations are indeed from an ICA generative model, then the minimization of the
pair-wise dependencies is sufficient to get the solution of the ICA task according
to the Darmois-Skitovich theorem \cite{comon94independent}. This is not the case for
the ISA task, however. There are ISA tasks, where the estimation of pair-wise
dependencies is insufficient for recovering the hidden subspaces
\cite{poczos05geodesic}. Nonetheless, such algorithms seem to work nicely in many
practical cases.

A further complication arises if the $d_e^m$ dimensions of subspaces $\b{e}^m$ are not known. Then the dimension of the
entropy estimation becomes uncertain.  Methods that try to apply pair-wise dependencies were proposed to this task. One
can find a block-diagonalization method in \cite{theis06towards}, whereas \cite{bach03finding} makes use of kernel
estimations of the mutual information.

Here we shall assume that the separation theorem is satisfied. We shall apply ICA
preprocessing. This step will be followed by the estimation of the pair-wise mutual
information of the ICA coordinates. These quantities will be considered as the
weights of a weighted graph, the vertices of the graph being the ICA coordinates. We
shall search for clusters of this graph. In our numerical studies, we make use of
Kernel Canonical Correlation Analysis \cite{bach03beyond} for the MI estimation. A
variant of the Ncut algorithm \cite{yu03multiclass} is applied for clustering. As a
result, the mutual information within (between) cluster(s) becomes large (small).

The problem is that this ISA method requires i.i.d. hidden sources. Below, we show how to generalize the ISA task to
more realistic sources. Finally, we solve this more general problem when the dimension of the subspaces is not known.

\section{ISA Generalizations}

We need the following notations: Let $z$ stand for the time-shift operation, that is
\mbox{$(z\b{v})(t):=\b{v}(t-1)$}. The N order polynomials of $D_1\times D_2$
matrices are denoted as \mbox{$\R[z]_N^{D_1\times
D_2}:=\{\b{F}[z]=\sum_{n=0}^N\b{F}_nz^n, \b{F}_n\in\R^{D_1\times D_2})\}$}. Let
$\nabla^r[z]:=(\b{I}-\b{I}z)^r$ denote the operator of the $r^{th}$ order
difference, where $\b{I}$ is the identity matrix,  $r\geq 0$, $r\in\mathbb{Z}$.

Now, we are to estimate unknown components  $\b{e}^m$ from
observed signals $\b{x}$. We always assume that $\b{e}$ takes the
form like in \eqref{eq:e_concat} and that $\b{A} \in \R^{D_x
\times D_s}$ is of full column rank.

\begin{enumerate}
    \item
        AR-IPA: The AR generalization of the ISA task is defined by the following equations:
        $\b{x}=\b{A}\b{s}$, where $\b{s}$ is an AR(p) process i.e, $\b{P}[z]\b{s}=\b{Q}\b{e}$,
        \mbox{$\b{Q} \in \R^{D_s \times D_e}$}, and $\b{P}[z]:=\b{I}_{D_s}-\sum_{i=1}^{p}\b{P}_{i}z^{i}\in \R[z]_{p}^{D_s\times
        D_s}$. We assume that $\b{P}[z]$ is stable, that is $\det(\b{P}[z] \neq 0)$, for all $z \in
        \C$, $|z| \leq 1$. For $d_e^m=1$ this task was investigated in \cite{hyvarinen98stochastic}.
         Case $d_e^m>1$ is treated in \cite{poczos05independent}. The special case of $p=0$ is the ISA task.
    \item
        MA-IPA or Blind Subspace Deconvolution (BSSD) task: The ISA task is generalized to blind deconvolution
        task (moving average task, MA(q)) as follows: $\b{x}=\b{Q}[z]\b{e}$, where
        $\b{Q}[z]=\sum_{j=0}^{q}\b{Q}_j z^j \in \R[z]_{q}^{D_x\times D_e}$.
    \item
        ARMA-IPA task:
         The two tasks above can be merged into a model, where the hidden $\b{s}$ is ARMA(p,q):
        $\b{x}=\b{A}\b{s}$, $\b{P}[z]\b{s}=\b{Q}[z]\b{e}$. Here $\b{P}[z]\in \R[z]_{p}^{D_s\times D_s}$,
        $\b{Q}[z]\in \R[z]_{q}^{D_s\times D_e}$. We assumed that $\b{P}[z]$ is stable. Thus the ARMA process is stationary.
    \item
        ARIMA-IPA task: In practice, hidden processes $\b{s}$ may be non-stationary. ARMA processes can be generalized
        to the non-stationary case. This generalization is called integrated ARMA, or ARIMA(p,r,q). The assumption
        here is that the $r^{th}$ difference of the process is an ARMA process. The corresponding IPA task is then
\begin{align}
\b{x}=\b{A}\b{s}, \text{ where }\b{P}[z]\nabla^r[z]\b{s}=\b{Q}[z]\b{e}.
\label{eq:ARIMA}
\end{align}
\end{enumerate}

\section{Reduction of ARIMA-IPA to ISA}\label{sec:ARIMA-IPA2ISA}

We show how to solve the above tasks by means of ISA algorithms. We treat the ARIMA
task. Others are special cases of this one. In what follows, we assume that: (i)
$\b{P}[z]$ is stable, (ii) the mixing matrix $\b{A}$ is of full column rank, and
(iii) $\b{Q}[z]$ has left inverse. In other words, there exists a polynomial matrix
\mbox{$\b{W}[z]\in\R[z]^{D_e\times D_s}$} such that $\b{W}[z]\b{Q}[z]=\b{I}_{D_e}$
\footnote{One can show for $D_s>D_e$ that under mild conditions $\b{Q}[z]$-has an
inverse with probability 1 \cite{rajagopal03multivariate}; e.g., when the matrix
$[\b{Q}_0,\ldots,\b{Q}_q]$ is drawn from a continuous distribution.}.

The route of the solution is elaborated here. Let us note that differentiating the
observation $\b{x}$ of the ARIMA-IPA task in Eq.~\eqref{eq:ARIMA} in $r^{th}$ order,
and making use of the relation \mbox{$z\b{x}=\b{A}(z\b{s})$}, the following holds:
\begin{align}
\nabla^r[z]\b{x}=\b{A}\left(\nabla^r[z]\b{s}\right), \text{ and }
\b{P}[z]\left(\nabla^r[z]\b{s}\right)&=\b{Q}[z]\b{e}.\label{eq:obs2}
\end{align}
That is taking $\nabla^r[z]\b{x}$ as observations, one ends up with an ARMA-IPA
task. Assume that $D_x>D_e$ (undercomplete case). We call this task uARMA-IPA. Now
we show how to transform the uARMA-IPA task to ISA. The method is similar to that of
\cite{gorokhov99blind} where it was applied for BSD.

\begin{theorem2}
If the above assumptions are fulfilled then in the uARMA-IPA task, observation process $\b{x}(t)$ is autoregressive and
its innovation $\tilde{\b{x}}(t):=\b{x}(t)-E[\b{x}(t)|\b{x}(t-1),\b{x}(t-2),\ldots]=\b{A}\b{Q}_0\b{e}(t)$, where
$E[\cdot| \cdot]$ denotes the conditional expectation value. Consequently, there is a polynomial matrix
\mbox{$\b{W}_{\mathrm{AR}}[z]\in\R[z]^{D_x\times D_x}$} such that $\b{W}_{\mathrm{AR}}[z]\b{x}=\b{A}\b{Q}_0\b{e}$.
\end{theorem2}

\begin{proof2}
Steps of the proof:
\begin{enumerate}
    \item
    In the uARMA-IPA task the following equations hold:
    \begin{align}
    \b{P}[z]\b{s}&=\b{Q}[z]\b{e}, \label{eq:hidden1b}\\
    \b{x}&=\b{A}\b{s},
    \end{align}
or equivalently
\begin{align}
    \b{s}(t)&=\sum_{i=1}^{p}\b{P}_i\b{s}(t-i)+\sum_{j=0}^{q}\b{Q}_j\b{e}(t-j),\label{eq:hidden1}\\
    \b{x}(t)&=\b{A}\b{s}(t). \label{eq:obs1}
\end{align}
        Non-degenerate linear transformation of an ARMA process is also ARMA. Thus, observation process $\b{x}$ is an ARMA process.
        Formally: Substituting $\b{s}(t)$ of Eq.~\eqref{eq:hidden1} into Eq.~\eqref{eq:obs1} and then using the pseudoinverse of matrix $\b{A}$
        and expression $\b{s}(t)=\b{A}^{-1}\b{x}(t)$ that follows from Eq.~\eqref{eq:obs1}, we have
        \begin{equation}
            \b{x}(t)=\sum_{i=1}^{p}\b{A}\b{P}_i\b{A}^{-1}\b{x}(t-i)+\sum_{j=0}^{q}\b{A}\b{Q}_j\b{e}(t-j).\label{eq:x-ARMA}
        \end{equation}
        Process $\b{e}(t)$ is i.i.d, so the process $\b{x}(t)$ is ARMA.
    \item
        We assumed that $\b{Q}[z]$ has left inverse and thus $\b{e}$ of Eq.~\eqref{eq:hidden1b} can be expressed from
        $\b{s}$ via multiplication with a polynomial matrix. One says that $\b{e}$ derives from $\b{s}$ by causal FIR filtering. The same holds for $\b{x}$ because of Eq.~\eqref{eq:obs1}:
        \begin{equation}
            \b{e}=\b{P}^{'}[z]\b{s}=\b{P}^{'}[z](\b{A}^{-1}\b{x})=(\b{P}^{'}[z]\b{A}^{-1})\b{x}=:\b{P}^{''}[z]\b{x},\label{eq:e-from-x-finite-past}
        \end{equation}
         \mbox{where $\b{P}^{'}[z]:=\b{Q}^{-1}[z]\b{P}[z]=\sum_{n=0}^N\b{P}^{'}_nz^{-n}\in\R[z]_N^{D_e\times D_s}$, $\b{P}^{''}[z]\in\R[z]_N^{D_e\times D_x}$} and $N$ denotes the degree of the polynomials.
    \item
        The first term of the r.h.s. of the observation $\b{x}$ in Eq.~\eqref{eq:x-ARMA} is a linear expression of a \emph{finite} history of
        $\b{x}$.  Equation~\eqref{eq:e-from-x-finite-past} implies, that the second term, except $\b{A}\b{Q}_0\b{e}(t)$, also belongs to the linear hull of the finite history of $\b{x}$.
        Formally:
        \begin{align}
            \b{x}(t)&=\b{A}\b{Q}_0\b{e}(t) + \sum_{i=1}^{p}\b{A}\b{P}_i\b{A}^{-1}\b{x}(t-i) +
            \sum_{j=1}^{q}\b{A}\b{Q}_j(\b{P}^{''}[z]\b{x})(t-j)\\
            &\in \b{A}\b{Q}_0\b{e}(t) + \left<\b{x}(t-1),\b{x}(t-2),\ldots,\b{x}(t-\max(p,q+N))\right>.
       \end{align}
    \item
        $\b{e}(t)$ is independent of $\left<\b{x}(t-1),\b{x}(t-2),\ldots,\b{x}(t-\max(p,q+N))\right>$. Consequently, observation process $\b{x}(t)$ is autoregressive with innovation
        $\b{A}\b{Q}_0\b{e}(t)$.
    \end{enumerate}
\end{proof2}

Thus, AR fit of $\b{x}(t)$ can be used for the estimation of
$\b{A}\b{Q}_0\b{e}(t)$. This innovation corresponds to the
observation of an undercomplete ISA model
($D_x>D_e$)\footnote{Assumptions made for $\b{Q}[z]$ and $\b{A}$
in the uARMA-IPA task implies that $\b{A}\b{Q}_0$ is of full
column rank and thus the resulting ISA task is well defined.},
which can be reduced to a complete ISA ($D_x=D_e$) using PCA.
Finally, the solution can be finished by any ISA procedure. The
reduction procedure implies that hidden components $\b{e}^m$ can
be recovered only up to the ambiguities of the ISA task:
components of (identical dimensions) can be recovered only up to
permutations. Within each subspaces, unambiguity is warranted only
up to orthogonal transformations.

The steps of our algorithm are summarized in Table
\ref{tab:pseudocode}.

\begin{table}
  \centering
  \caption{Pseudocode of the undercomplete ARIMA-IPA algorithm} \label{tab:pseudocode}
  \begin{tabular}{|l|}
        \hline
        \textbf{Input of the algorithm}\\
        \verb|   |Observation: $\{\mathbf{x}(t)\}_{t=1,\ldots,T}$\\
        \textbf{Optimization}\\
        \verb|   |\textbf{Differentiating}: for observation $\b{x}$ calculate $\b{x}^*=\nabla^r[z]\b{x}$ \\
        \verb|   |\textbf{AR fit}: for $\b{x}^*$ estimate $\hat{\b{W}}_{\text{AR}}[z]$\\
        \verb|   |\textbf{Estimate innovation}: $\tilde{\b{x}}=\hat{\b{W}}_{\text{AR}}[z]\b{x}^*$\\
        \verb|   |\textbf{Reduce uISA to ISA and whiten}: $\tilde{\b{x}}^{'}=\hat{\b{W}}_{\text{PCA}}\tilde{\b{x}}$\\
        \verb|   |\textbf{Apply ICA for $\tilde{\b{x}}^{'}$}: $\b{e}^{*}=\hat{\b{W}}_{\text{ICA}}\tilde{\b{x}}^{'}$\\
        \verb|   |\textbf{Estimate pairwise dependency} e.g., as in \cite{bach03finding} on $\b{e}^{*}$\\
        \verb|   |\textbf{Cluster $\b{e}^{*}$ by Ncut}: the permutation matrix is $\mathcal{P}$ \\
        \textbf{Estimation}\\
        \verb|    |$\hat{\b{W}}_{\text{ARIMA-IPA}}[z]=\mathcal{P}\hat{\b{W}}_{\text{ICA}}\hat{\b{W}}_{\text{PCA}}\hat{\b{W}}_{\text{AR}}[z]\nabla^r[z]$\\
        \verb|               |$\hat{\b{e}}=\hat{\b{W}}_{\text{ARIMA-IPA}}[z]\b{x}$\\
        \hline
 \end{tabular}
\end{table}

\section{Results}

In this section we demonstrate the theoretical results by numerical simulations.

\subsection{ARIMA Processes}
\begin{figure}[h!]
\centering%
\subfloat[][]{\includegraphics[width=1cm]{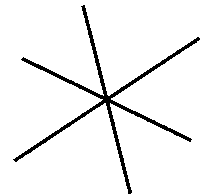}
\includegraphics[width=.9cm]{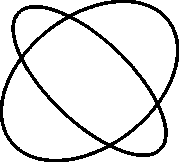}\hspace*{0.08cm}
\includegraphics[width=1cm]{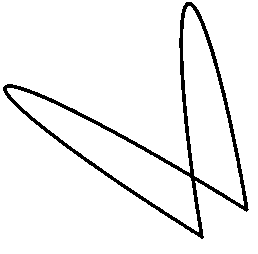}\label{fig:datbase3d}}\hspace*{0.5cm}
\subfloat[][]{\includegraphics[width=4.85cm]{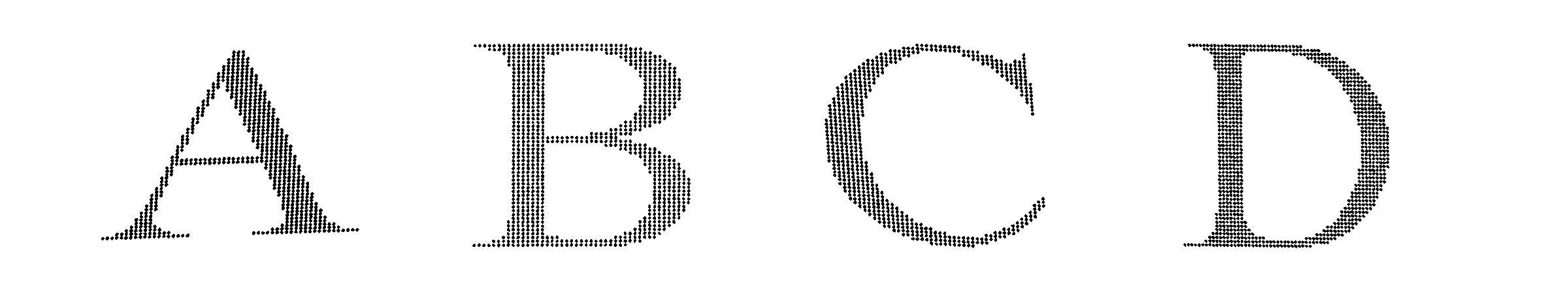}\label{fig:datbase2d}}\hspace*{0.65cm}
\subfloat[][]{\includegraphics[width=1cm]{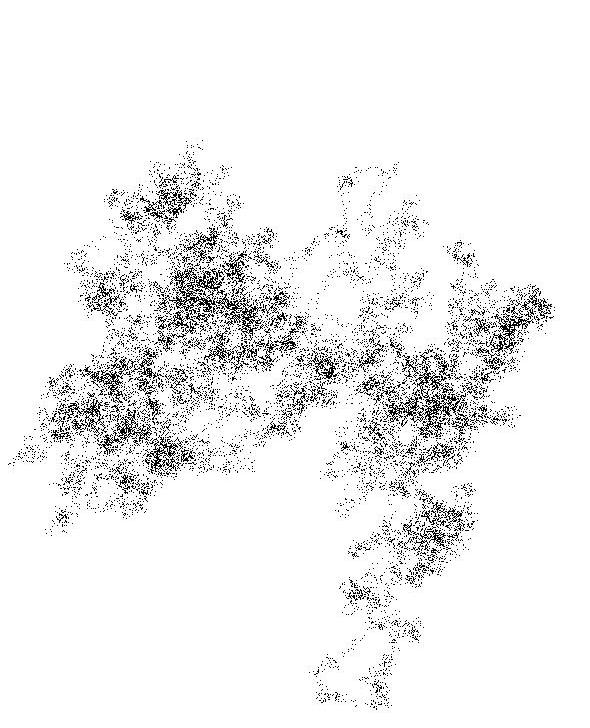}\label{fig:obs3d}}\hfill
\subfloat[][]{\includegraphics[width=1cm]{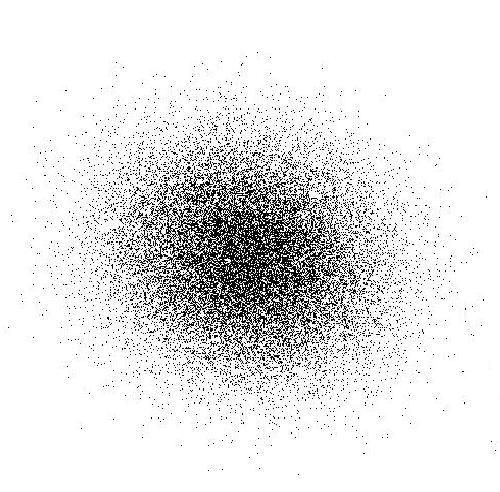}\label{fig:diff_obs3d}}\\
\hspace*{-0.06cm}
\subfloat[][]{\includegraphics[width=3.2cm]{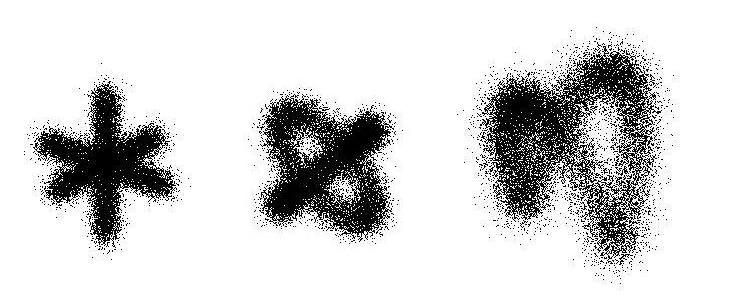}\label{fig:est3d}}\hspace*{0.73cm}
\subfloat[][]{\includegraphics[width=4.5cm]{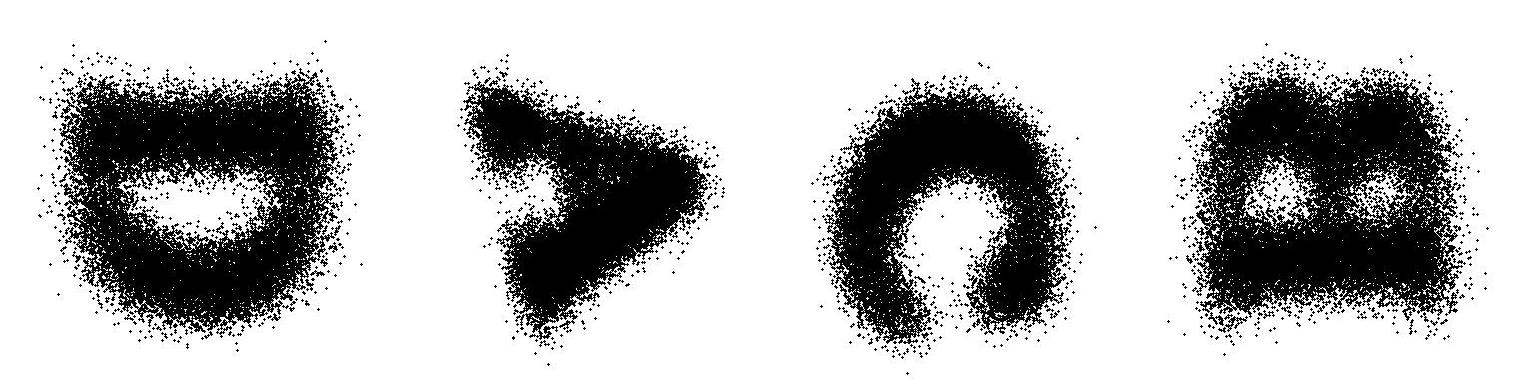}\label{fig:est2d}}\hspace*{1.19cm}
\subfloat[][]{\includegraphics[width=2cm]{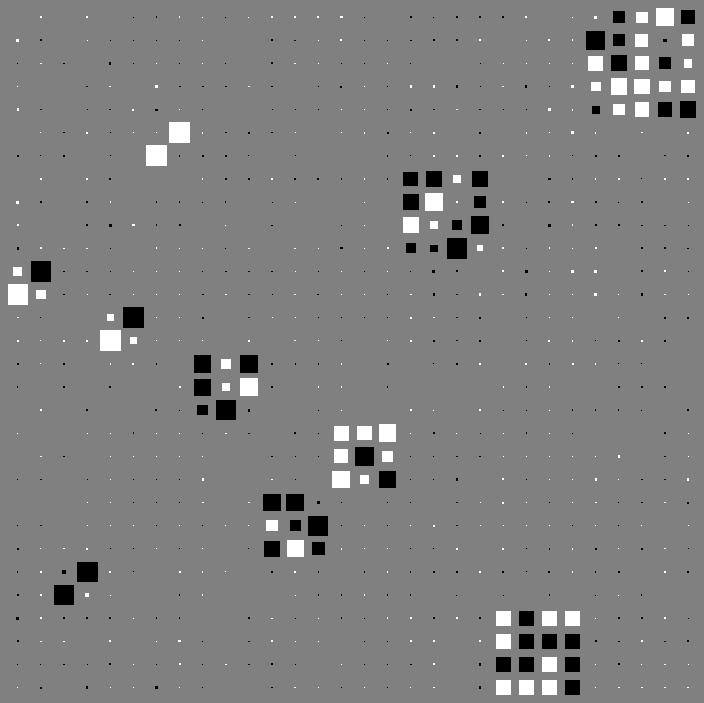}\label{fig:hinton}} \caption[]{
(a-b) components of the database. (a): 3 pieces of 3D geometrical forms, (b): 4
pieces of 2D letters. Hidden sources are uniformly distributed variables on these
objects. (c): typical 3D projection of the observation. (d): typical 3D projection
of the $r^{th}$-order difference of the observation, (e): estimated 2D components,
(f): estimated 3D components, (g): Hinton diagram of $\b{G}$, which -- in case of
perfect estimation -- becomes a block permutation matrix.}
\end{figure}

We created a database for the demonstration: Hidden sources $\b{e}^m$ are 4 pieces
of 2D, 3 pieces of 3D, 2 pieces of 4D and 1 piece of 5D stochastic variables, i.e.,
$M=10$. These stochastic variables are independent, but the coordinates of each
stochastic variable $\b{e}^m$ depend on each other. They form a 30 dimensional space
together $(D_e=30)$. For the sake of illustration, 3D (2D) sources emit random
samples of uniform distributions defined on different 3D geometrical forms (letters
of the alphabet). The distributions are depicted in Fig.~\ref{fig:datbase3d}
(Fig.~\ref{fig:datbase2d}). 30,000 samples were drawn from the sources and they were
used to drive an ARIMA(2,1,6) process defined by \eqref{eq:ARIMA}. Matrix $\b{A} \in
\R^{60 \times 60}$ was randomly generated and orthogonal. We also generated
polynomial $\b{Q}[z]\in \R[z]_5^{60 \times 30}$ and stable polynomial $\b{P}[z] \in
\R[z]_1^{60 \times 60}$ randomly. The visualization of the 60 dimensional process is
hard to illustrate: a typical 3D projection is shown in Fig.~\ref{fig:obs3d}. The
task is to estimate original sources $\b{e}^m$ using these non-stationary
observations. $r^{th}$-order differencing of the observed ARIMA process gives rise
to an ARMA process. Typical 3D projection of this ARMA process is shown
Fig.~\ref{fig:diff_obs3d}. Now, one can execute the other steps of
Table~\ref{tab:pseudocode} and these steps provide the estimations of the hidden
components $\hat{\b{e}}^m$. Estimations of the 3D (2D) components are provided in
Fig.~\ref{fig:est3d} (Fig.~\ref{fig:est2d}). In the ideal case, the product of
matrix $\b{AQ}_0$ and the matrices provided by PCA and ISA, i.e.,
$\b{G}:=(\mathcal{P}\hat{\b{W}}_{\text{ICA}}\hat{\b{W}}_{\text{PCA}})\b{AQ}_0\in\R^{D_e\times
D_e}$ is a block permutation matrix made of $d^m_e\times d^m_e$ blocks. This is
shown in Fig.~\ref{fig:hinton}.

\vspace*{-.5cm}

\subsection{Facial Components}

\begin{figure}[h!]%
\centering%
\subfloat[][]{
\includegraphics[width=5cm]{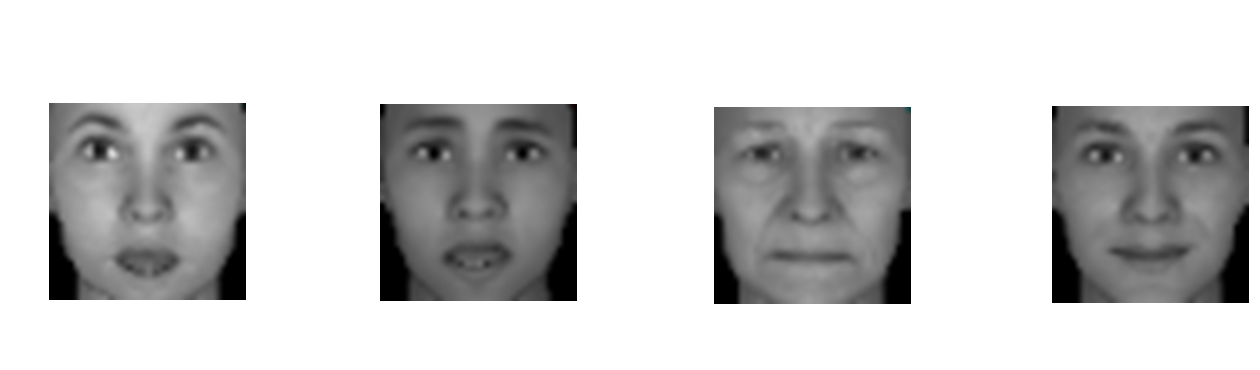}\hfill \label{fig:facedatabase}} \\
\subfloat[][]{\includegraphics[width=12cm]{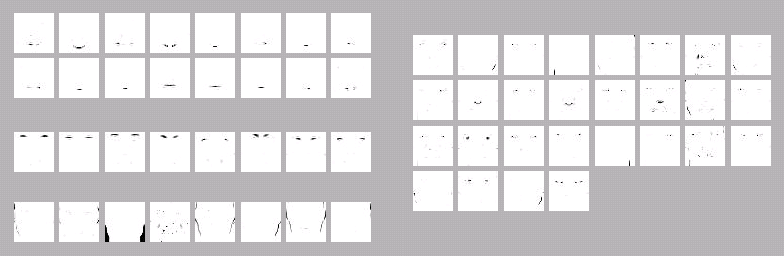} \label{fig:kernelmi}
\hfill} \caption[]{ (a) Samples from the database. (b) Four subspaces of the
components. Components in distinct groups correspond mostly to
mouth, eye brushes, facial profiles, and eyes respectively.}%
\label{fig:database}%
\end{figure}

We have generated another database using the
FaceGen\footnote{\url{http://www.facegen.com/modeller.htm}} animation software. In
our database we had 800 different front view faces with the 6 basic facial
expressions. We had thus 4,800 images in total. All images were sized to $40 \times
40 $ pixel. Figure~\ref{fig:facedatabase} shows samples of the database. A large
$\b{X} \in \R^{4800 \times 1600}$ matrix was compiled; rows of this matrix were 1600
dimensional vectors formed by the pixel values of the individual images.  The
\emph{columns} of this matrix were considered as mixed signals. This treatment
replicates the experiments in \cite{bartlett02face}: Bartlett et al., have shown
that in such cases, undercomplete ICA finds components resembling to what humans
consider facial components. We were interested in seeing the components grouped by
undercomplete ISA algorithm. The observed 4800 dimensional signals were compressed
by PCA to $60$ dimensions and we searched for 4 pieces of ISA subspaces using the
algorithm detailed in Table~\ref{tab:pseudocode}. The 4 subspaces that our algorithm
found are shown in Fig.~\ref{fig:kernelmi}. As it can be seen, the 4 subspaces
embrace facial components which correspond mostly to mouth, eye brushes, facial
profiles, and eyes, respectively.

\section{Conclusions}

We have extended the ISA task to problems where the hidden components can be AR, MA,
ARMA, or ARIMA processes. We showed an algorithm that can identify the hidden
subspaces under certain conditions. The algorithm does not require previous
knowledge about the dimensions of the subspaces. The working  of the algorithm was
demonstrated on an artificially generated ARIMA process, as well as on a database of
facial expressions.

%\section*{Acknowledgment}

\bibliographystyle{unsrt}
%\bibliography{ICA-2007}

\end{document}